\documentclass[12pt,reqno]{amsart}
\usepackage{geometry}                		
\usepackage{amsmath}
  \usepackage{paralist}
  \usepackage{graphics} 
  \usepackage{epsfig} 
\usepackage{graphicx}  \usepackage{epstopdf}
 \usepackage[colorlinks=true]{hyperref}
\hypersetup{urlcolor=blue, citecolor=red}

\numberwithin{equation}{section}
\allowdisplaybreaks[4]
\newcommand{\xqedhere}[2]{%
	\rlap{\hbox to#1{\hfil\llap{\ensuremath{#2}}}}}

\newtheorem{theorem}{Theorem}[section]

\newtheorem{lemma}[theorem]{Lemma}
\newtheorem{proposition}{Proposition}

\theoremstyle{definition}
\newtheorem{definition}[theorem]{Definition}
\newtheorem{remark}{Remark}

\title[On the invariant region for compressible Euler equations] 
      {On the invariant region for compressible Euler equations with a general equation of state}

\author[Hailiang Liu and Ferdinand Thein]{}

\dedicatory{This article has been published in a revised form in \emph{Communications on Pure and Applied Analysis} \url{doi: 10.3934/cpaa.2021084}.
This version is free to download for private research and study only. Not for redistribution, re-sale or use in derivative works.\\
~\newline
Dedicated to Professor Shuxing Chen on the occasion of his 80th birthday}



 \email{hliu@iastate.edu}
 \email{ferdinand.thein@ovgu.de}

\thanks{Hailiang Liu was partially supported by the National Science Foundation under Grant DMS1812666}

\thanks{$^*$ Corresponding author: Hailiang Liu}

\begin{document}

\maketitle

\centerline{\scshape Hailiang Liu$^*$}
\medskip
{\footnotesize
 \centerline{Iowa State University, Mathematics Department}
   \centerline{Ames, IA 50011, USA}
} 

\medskip

\centerline{\scshape Ferdinand Thein}
\medskip
{\footnotesize
 \centerline{Otto-von-Guericke-Universit\"at, Universit\"atsplatz 2}
   \centerline{Magdeburg, 39106, Germany}
}

\bigskip


\begin{abstract}
    The state space for solutions of the compressible Euler equations with a general equation of state is examined.
    An arbitrary equation of state is allowed, subject only to the physical requirements of thermodynamics.
    An invariant region of the resulting Euler system is identified and the convexity property of this region is justified by using only very minimal thermodynamical assumptions.
    Finally, we show how an invariant-region-preserving (IRP) limiter can be constructed for use in high order finite-volume type schemes
    to solve the compressible Euler equations with a general constitutive relation.
\end{abstract}

\section{Introduction}
The dynamical evolution of a fluid is determined by the principles of conservation of mass, momentum, and energy. 
Material properties strongly influence the structure and dynamics of waves in any continuum-mechanical system, as shown by Menikoff and Plohr in their instructive and comprehensive work \cite{MP89}.
There the Riemann problem for the compressible Euler equations is discussed in detail.
Constitutive relations are understood as the most fundamental level for fluids, where they are embodied in an equilibrium equation of state (EOS)
and the physical principles of thermodynamics impose stringent constraints on the equation of state, cf.\ Landau and Lifshitz \cite{Landau1987}.

In this work we investigate the invariant region for solutions to the Euler equations with a general EOS.
This system arises in fluid dynamics and is the paradigm for systems of hyperbolic conservation laws.
It provides the motivation for many of the central ideas in the analysis of quasilinear hyperbolic partial differential equations. The Euler equations we consider here are given by
\begin{subequations}\label{eu}
    \begin{align}
        & \partial_t \rho +\nabla_x\cdot(\rho \mathbf{u})=0,\\
        & \partial_t(\rho \mathbf{u}) +\nabla_x \cdot(\rho \mathbf{u}\otimes \mathbf{u}) +\nabla_x P=0,\\
        & \partial_t E +\nabla_x\cdot((E+P)\mathbf{u})=0, 
    \end{align}
\end{subequations}
where $\rho$ denotes the mass density, $\mathbf{u}$ the fluid velocity, $E$ the specific total energy, and $P$ the pressure.
The pressure must be specified through a constitutive relation, the EOS, which characterizes the fluid material.
The qualitative character of the solutions of the conservation laws depends crucially on the thermodynamic properties of the material.
For such a general system as we consider it here, a complete mathematical theory is still unavailable.
A natural choice would be to consider solutions confined in a convex domain in state space,
which is invariant either provable to hold in time along dynamic equations or known from physical considerations.
An invariant region to (\ref{eu}) is an open set in phase space $\mathbb{R}^d$ such that if initial data lie in this set, then the solution will remain in this set.
It was proved by Hoff \cite{Hoff85} that an invariant region for one dimensional hyperbolic conservation laws must be convex.
The invariant region is important also for numerical algorithms to solve these equations.
However, there are no rigorous convergence results to solutions for any numerical schemes approximating multi-dimensional systems of conservation laws.
Nevertheless, if the numerical solution lies in some important invariant region, this would be an indication that the problem is well-posed or the solution is physically meaningful.
Some first order methods such as the Lax-Friedrich method and the Godunov method are known to satisfy the invariant-region-preserving (IRP) property \cite{Fr95, Fr01, Ta86}.
In the context of continuous finite elements,
the IRP property has been studied by Guermond and Popov \cite{GP16} using the first order approximation to solve general hyperbolic conservation law systems.

However, it is a rather difficult task to preserve an invariant region by a high order numerical method unless some nonlinear limiter is frequently imposed;
cf.\, second-order limitation techniques are used in \cite{KP94} to enforce the minimum entropy principle so that oscillations in numerical solutions may be damped.
Indeed, recent efforts have been using limiting techniques to construct high order schemes, including maximum-principle-preserving schemes for scalar conservation laws (see \cite{ZhangShu2010a})
and positivity-preserving schemes for hyperbolic systems including the compressible Euler equations (see e.g. \cite{PerthameShu1996, ZhangShu2010b, ZhangXiaShu2012}).
The work by Zhang and Shu in \cite{ZhangShu2012} introduced a limiter to preserve the minimum-entropy-principle.

The emphasis in \cite{JL18a, JL18b, JL19} is on the notion of invariant regions and how to numerically preserve an invariant region by high order discontinuous Galerkin schemes.
We observe that the ideas for the high order IRP schemes studied in \cite{JL18b} can be extended to general hyperbolic conservation law systems as long as it features a convex invariant region.
In this work we show how to identify a convex invariant region for the Euler system with a general EOS
and present a simple IRP limiter which when coupled with any high order finite volume type methods can pull numerical solutions back to the invariant region.
It is hoped that this information will be of service to the designers of numerical approximations of this important class of equations.

The remainder of the paper is organized as follows. In the next section, we review special models of the Euler equations and present the known invariant region for each of them.
Further we discuss the basic structure for the full system of Euler equations, and present the main result to be proved in Section 3 and 4.
In Section 3, we show a region defined by the entropy to be invariant. In Section 4, we provide a justification for a sub-level set of the internal energy.
Section 5 is devoted to the IRP limiter and discussion on high order IRP schemes. Finally, a specific example of the general Euler system is given in Section 6.

\section{Equation of state}
To obtain a closed system it is necessary to prescribe an equation of state, which relates the pressure $P$ to the density and the total energy.
It is provided by thermodynamics. The following three distinct situations are important:

\subsection{Pressureless gases}
The pressure $P$ vanishes and so the total energy reduces to just the kinetic energy: $E=\frac{1}{2}\rho |\mathbf{u}|^2$. The equations (\ref{eu}) take the form
\begin{align*}
    & \partial_t \rho +\nabla_x\cdot(\rho \mathbf{u})=0,\\
    & \partial_t(\rho \mathbf{u}) +\nabla_x \cdot(\rho \mathbf{u}\otimes \mathbf{u})=0,
\end{align*}
and the energy equation in (\ref{eu}) follows formally from the continuity and momentum equations.
The pressureless Euler system has been proposed as a simple model describing the formation of galaxies in the early stage of the universe.
There is a maximum principle for the velocity field $\mathbf{u}$, though its gradient can become discontinuous -- shock formation.
If fluid elements meet at the same location, then they stick together to form larger compounds and so $\rho(x, t)$ can have singular parts (in particular, Dirac measures).
Consequently the system must be understood in the sense of distributions, and admits an invariant region of the form
$$
\Sigma=\{(\rho, \mathbf{m})|\quad \rho>0, \quad -\alpha_i \rho \leq m_i \leq \alpha_i \rho, i=1,\dots, d\}
$$
for some $\alpha_i>0$ with $\alpha_i=\sup_{x} |u_i(x, 0)|$, and $\mathbf{m}=\rho \mathbf{u}$.

\subsection{Isentropic gases}
In this regime, the thermodynamical entropy of the fluid is assumed to be constant in space and time. Consequently, the pressure is a function of the density only. For polytropic gases,
$$
P(\rho)=U'(\rho)\rho -U(\rho), \quad \rho \geq 0,
$$
so that $\rho e = U(\rho)$ is convex in $\rho$ and $E=\frac{1}{2}\rho|\mathbf{u}|^2 +U(\rho)$.
As in the pressureless case, the energy equation in (\ref{eu}) follows formally from the continuity and the momentum equation.
For solutions with shocks, the continuity and the momentum equation must be considered in the sense of distributions, and the energy equation does no longer follow automatically.
A physically reasonable relaxation is to assume that no energy can be created by the fluid: the energy equality in (\ref{eu}) must be replaced by the inequality
$$
\partial_t \left(\frac{1}{2} \rho |\mathbf{u}|^2 +U(\rho)\right) +\nabla\cdot \left(\left( \frac{1}{2}\rho |\mathbf{u}|^2 +U'(\rho)\rho \right) \mathbf{u} \right) \leq 0
$$
in distributional sense. If the above inequality is strict, this means physically that mechanical energy is transformed into heat.
A form of energy that is not accounted for by the model. For one dimensional flows, we find Riemann invariants of the form
$$
R_\pm =\frac{m}{\rho} \pm \int^\rho \sqrt{\frac{U''(r)}{r}}dr,
$$
so that for $\lambda_\pm= u \pm \sqrt{U''(\rho)\rho}$,
\begin{align*}
    & \partial_t R_- +\lambda_- \partial_x R_- =0, \\
    &  \partial_t R_+ +\lambda_+ \partial_x R_+ =0.
\end{align*}
The corresponding invariant region is given by
$$
\Sigma=\{ (\rho, m)^\top|\quad R_+ \leq c_2, \quad R_- \geq c_1\},
$$
where $c_i$ is determined by the given initial data.

\subsection{Full Euler equations}
We consider a polytropic gas with the adiabatic coefficient $\gamma_0 > 1$. Then the pressure is given in terms of the conserved quantities $\mathbf{w} :=(\rho, \mathbf{m}, E)^\top$ by the formula
$$
P= (\gamma_0-1)\left(E-\frac{|\mathbf{m}|^2}{2\rho}\right) = (\gamma_0-1)\left(E-\frac{1}{2}\rho |\mathbf{u}|^2\right) = (\gamma_0 - 1)\rho e.
$$
Density and pressure define the specific thermodynamical entropy, given as
$$
s:={\rm log} \left( \frac{P}{c\rho^{\gamma_0}}\right), \quad c=k(\gamma_0-1)>0
$$
in the case of polytropic gases. The internal energy density has the form
\[
  e= k\exp(s) \rho^{\gamma_0-1}.
\]
The specific entropy $s$ must be constant along characteristics:
$$
\partial_ t s +\mathbf{u} \cdot \nabla_x s=0.
$$
But since the solutions to the compressible Euler equations may become discontinuous in finite time, the physically reasonable relaxation is
that the specific entropy should be non-decreasing forward in time, which expresses the second law of thermodynamics. It follows that
$$
\inf s(x, t) \geq \inf s(x, 0)=:s_0, 
$$
where $s(x, 0)$ is the initial specific entropy. Hence the region
$$
\Sigma=\{\mathbf{w}\;|\; \rho > 0,\; P>0,\; s\geq s_0\}
$$
is invariant and convex (see e.g.,\ \cite{JL18b}).

\subsection{General EOS}
From now on, we consider a general constitutive relation derived from thermodynamics.
By the second law of thermodynamics we know that there exists a function $s$, called entropy, which is a twice differential function of form
$$
s=G(e, v),
$$
where $e$ is the specific internal energy and $v=1/\rho > 0$ is the specific volume. This formula means
$$
s(x, t)=G(e(x, t), v(x, t)),
$$
for $t>0$ and $x$ in the fluid region. The thermodynamic properties of a material are embodied in this relation.
This fundamental equation contains all the information we need to close the Euler system. As usual, we assume that $s$ has the following properties:
\begin{enumerate}[(I)]
    \item $G(e, v)$ is concave in $(e, v)$,
    \item $\partial_e G(e, v)>0$.
\end{enumerate}
As a result, $s=G(e, v)$ can be reformulated as 
$$
e=F(s, v),
$$
which for $G$ strict concave can be shown to have the following properties:
\begin{enumerate}[(I)]
    \item $F(s, v)$ is convex in $(s, v)$, 
    \item $\partial_s F(s, v)=\dfrac{1}{G_e}>0$.
\end{enumerate}
We can define the pressure $P$ and the temperature $\theta$ by 
$$
P=-F_v, \quad \theta = F_s,
$$
in accordance with the fundamental thermodynamic identity 
$$
de=-Pdv +\theta ds,
$$
i.e., the first law of thermodynamics (see e.g.,\ \cite{Ev}). Imposed by the  Clausius-Duhem inequality, the entropy inequality and the energy equation take the form 
$$
0\leq \rho(\theta -F_s) (\partial_t s +\mathbf{u}\cdot \nabla_x s) +(-P-F_v)\nabla_x \mathbf{u}, 
$$
which must hold for all admissible thermodynamic processes. The above choice is a valid candidate for dynamic equations and also consistent with classical formulas for a simple fluid.
Rewriting the above we also have
$$
ds=\frac{1}{\theta}de +\frac{P}{\theta}dv.
$$
Hence  $\theta=1/G_e$ and $P=\theta G_v$ and thus we have
\begin{align}
    P=\frac{G_v}{G_e}.\label{eqdef:pressure}
\end{align}
This formula for the pressure enables us to close the Euler system of equations.\\
For our purposes we further follow Menikoff and Plohr and define dimensionless quantities as in \cite{MP89}.
In the situation under consideration it is possible to describe the fundamental properties of the studied system using these three quantities.
\begin{definition}[Dimensionless quantities]\label{def:diml_quantities}
    Given a thermodynamic system described by the fundamental equation $e = F(s,v)$ the following dimensionless quantities can be defined:
    \begin{enumerate}[(i)]
        \item The \emph{adiabatic exponent} or synonymous the \emph{dimensionless sound speed}
        \begin{equation}
            \gamma = \frac{v}{P}\frac{\partial^2 F}{\partial v^2}.\label{adib_exp}
        \end{equation}
        \item The \emph{Gr\"uneisen coefficient}
        \begin{equation}
            \Gamma = -\frac{v}{\theta}\frac{\partial^2 F}{\partial s \partial v}.\label{grueneisen_coeff}
        \end{equation}
        \item The \emph{dimensionless specific heat}
        \begin{equation}
            g = \frac{Pv}{\theta^2}\frac{\partial^2 F}{\partial s^2}.\label{diml_spec_heat}
        \end{equation}
        \item The \emph{fundamental derivative}
        \begin{equation}
            \mathcal{G} = -\frac{1}{2}v\dfrac{\dfrac{\partial^3 F}{\partial v^3}}{\dfrac{\partial^2 F}{\partial v^2}}
            = -\frac{1}{2}v\dfrac{\dfrac{\partial^2 P}{\partial v^2}}{\dfrac{\partial P}{\partial v}}.\label{fundamental_deriv}
        \end{equation}
    \end{enumerate}
\end{definition}
These dimensionless quantities are very helpful in discussing the Riemann problem for the Euler equations. Further they have several useful interpretations which we will not present here in detail.
We therefore again highly recommend \cite{MP89} and just make some remarks.
\begin{remark}\label{rem:diml_quantities}
    ~
    \begin{enumerate}[(i)]
        \item Systems with a convex fundamental equation $F(s,v)$ are called thermodynamically stable.
        \item The convexity of $F(s,v)$ can be expressed in terms of $g, \gamma$ \text{and} $\Gamma$, i.e.
        \begin{equation}
            F(s,v)\;\text{is  convex}\quad\Leftrightarrow\quad g \geq 0, \gamma \geq 0\;\text{and}\;g\gamma \geq \Gamma^2.\label{diml:E_convex}
        \end{equation}
        \item The thermodynamic laws do not confine the sign of $\mathcal{G}$.
        \item If $\gamma > 0$ holds, the system of Euler equations is strictly hyperbolic.
        \item The characteristic families corresponding to the largest and smallest eigenvalue of the full Euler system are genuine nonlinear iff $\mathcal{G} \neq 0$ holds.
        In particular shock waves are compressive iff $\mathcal{G} > 0$.
        \item For the polytropic gas presented above we obtain (cf.\,\cite{MP89})
        \begin{equation*}
            \gamma = \gamma_0,\,g = \gamma_0 - 1,\,\Gamma = \gamma_0 - 1\;\text{and}\;\mathcal{G} = \frac{\gamma_0 + 1}{2}.
        \end{equation*}
        Thermodynamic stability requires $\gamma_0 \geq 1$ and thus using this EOS the system of Euler equations is strictly hyperbolic.
    \end{enumerate}
\end{remark}
For further useful remarks and insights we refer to \cite{MP89}. Assuming $\mathcal{G} > 0$ and convexity of $F(s,v)$ it is possible to prove the following useful inequality
\[
 \frac{Pv}{e} < 2\gamma.
\]
The proof relies on the convexity of the isentropes in the $P-v$-plane, see \cite{MP89} (p.\ 95, Lem.\ 4.3). In view of equation (\ref{eqdef:pressure}) this gives
\begin{equation}
    \frac{vG_v}{eG_e} < 2\gamma.\label{ineq:pve_gamma}
\end{equation}
From now on we assume the convexity of the fundamental equation, see (\ref{diml:E_convex}) and additionally we demand $\mathcal{G} > 0$.
For consequences when $\mathcal{G} \leq 0$ we refer to M\"uller and Voss \cite{MV2006} and again to \cite{MP89}.

We recall that for $2\times 2$ systems such as the isentropic Euler system, an invariant region can be described by two Riemann invariants \cite{Sm83}.
For the Euler equations (\ref{eu}) with a general EOS, we shall identify a region which is both convex and invariant.
The convexity may be determined by the structure conditions of the EOS, and the invariance property depends on the dynamic equations. The main result can be stated as follows.
\begin{theorem}
    Consider system (\ref{eu}) and let the specific internal energy be given by $e = F(s,v)$ with the properties (I) and (II) given above.
    The quantities $\gamma,\ g,\ \Gamma$ and $\mathcal{G}$ are given as in Definition \ref{def:diml_quantities}. Furthermore the following constraint holds
    \[
      \mathcal{G} > 0.
    \]
    Then there exists a convex invariant region for the system of Euler equations:
    \begin{equation}\label{ir}
        \Sigma=\{\mathbf{w}=(\rho, \mathbf{m}, E)\,|\quad \rho>0,\; R>0,\;  q<0\},
    \end{equation}
    where
    $$
    R=E- \frac{|\mathbf{m}|^2}{2\rho}, \quad q=\rho (\inf s_0 - s).
    $$
\end{theorem}
Throughout the following two sections we will prove this theorem.
\section{Entropy minimization}
Using $e_s=\theta $ and $e_\rho=P/\rho^2$, and the energy equation, one can derive
\begin{equation*}
    0 =\partial_t E+\nabla_x \cdot((E+P) \mathbf{u})
    =\theta(\partial_t (\rho s)+\nabla_x \cdot(\rho s \mathbf{u})),
\end{equation*}
leading to the transport equation
$$
\partial_t s +\mathbf{u} \cdot \nabla_x s=0,
$$
at least for solution without shocks. For weak solutions, we need to recall the celebrated Clausius-Duhem inequality
$$
\partial_t (\rho s)+\nabla_x \cdot(\rho s \mathbf{u}) \geq \frac{r \rho}{\theta}-\nabla_x \cdot \left(\frac{f}{\theta} \right),
$$
which, when both the heat supply $r$ and heat flux $f$ vanish, leads to 
$$
\partial_t s +\mathbf{u} \cdot \nabla_x s\geq 0.
$$
In either case, we have the following entropy minimization principle:
$$
{\rm inf} s(x, t) ={\rm inf} s(x, 0)=:s_0, \quad s(x, t) -s_0 \geq 0.
$$
This ensures that the region $\{\mathbf{w}\,|\; \rho>0,\; q\leq 0\}$ with $q=\rho(s_0-s)$  is invariant.
Moreover, we have the following.
\begin{proposition}
    Let $q=\rho(s_0-s)$, then 
    $q$ is convex in terms of $\mathbf{w}:=(\rho, m, E)$ with $m=|\mathbf{m}|$.  
\end{proposition}
\begin{proof}
    For simplicity of computation, we take $s_0=0$ so that $ q=-\rho G(e, v).$
    A direct calculation gives
    \begin{align*}
        q_\rho & =-G +v G_e  \xi  + v G_v,  \; \xi=E- vm^2\\
        q_m & =mv G_e, \\
        q_E & =-G_e
    \end{align*}
    and second order derivatives 
    \begin{align*}
        q_{\rho \rho} &= -v^3 (G_{vv} +G_{ee} \xi^2 +2\xi  G_{ev} -m^2 G_e),\\
        q_{\rho m} & =-mv^2 G_e - G_{ee}v^3 m \xi -mv^3 G_{ev},\\
        q_{E\rho} & = v^2 \xi G_{ee} +v^2 G_{ev}, \;\\
        q_{mm} & = vG_e -v^3 m^2 G_{ee},\\
        q_{mE} &=v^2 m G_{ee},\\
        q_{EE} & = -vG_{ee},  \\
        q_{E\rho} & = v^2 \xi G_{ee} +v^2 G_{ev}.
    \end{align*}
    Note that from concavity of $s=G(e, v)$ we have
    \begin{equation}\label{ss}
        G_{vv}<0, \quad G_{ee}<0, \quad G_{vv} G_{ee} >G_{ev}^2.
    \end{equation}
    We shall only use this basic assumption to show $q$ is convex in $\mathbf{w}$. In order to do so, we only need to show
    $$
    q_{\rho\rho}>0, \quad A=q_{\rho\rho} q_{mm} -q_{m \rho}^2>0 \; \text{and}\;  B={\rm det} (D^2 q)>0.
    $$
    First from $q_{\rho \rho}$ we have
    \begin{align*}
        q_{\rho \rho}  & \geq -v^3 \left( G_{vv} +G_{ee}\xi^2 +2|\xi|\sqrt{G_{vv}G_{ee}} -m^2 G_e \right)\\
        &  = v^3m^2 G_e +v^3(\sqrt{|G_{vv}|}+|\xi|\sqrt{|G_{ee}|})^2 
         >0.
    \end{align*}
    Next we estimate $A$ as follows.
    \begin{align*}
        A & =q_{\rho \rho} q_{mm} -q_{\rho m}^2 \\ 
        & =- v^4 (G_e \!-\!v^2 m^2 G_{ee}) (G_{vv} \!+\!G_{ee}\xi^2 \!+\!2\xi G_{ev} \!-\!m^2 G_e)
       \!-\!m^2 v^4 (G_e \!+\!v\xi G_{ee} \!+\!vG_{ev})^2.
    \end{align*}
    Regrouping terms so that
    \begin{align*}
        -\rho^4 A & = G_e (G_{vv} +G_{ee}\xi^2 +2\xi G_{ev} -m^2 G_e +m^2 G_e  +2vm^2 (\xi G_{ee} +G_{ev}))\\
        & \quad -m^2 v^2 G_{ee} ( G_{vv} +G_{ee}\xi^2 +2\xi G_{ev} -m^2 G_e)\\
        & \quad +v^2 m^2 (\xi^2 G^2_{ee} +G^2_{ev} +2\xi G_{ee}G_{ev}) \\
        & =  G_e (E^2 G_{ee} +2E G_{ev}+G_{vv}) + v^2m^2 (G_{ev}^2 -G_{ee}G_{vv}) \\
        & =- G_e(E \sqrt{|G_{ee}|} +\sqrt{|G_{vv}|})^2 -v^2m^2(G_{ee}G_{vv}-G_{ev}^2)
         <0.
    \end{align*}
    Next we turn to estimate $B$. Note that
    \begin{align*}
        B -q_{EE} A & = 2q_{\rho E} q_{\rho m}q_{mE} -q_{mm}q_{\rho E}^2 -q_{\rho \rho }q_{mE}^2\\
        & = -2m^2 v^6 (G_{ev}+\xi G_{ee}) \cdot G_{ee} \cdot (G_e +\xi vG_{ee}+vG_{ev}) \\
        & \quad - v^6 (\rho G_e - v m^2 G_{ee}) \cdot (G_{ev}+\xi G_{ee})^2 \\
        & \quad + m^2 v^7 G_{ee}^2 (G_{vv} +\xi^2 G_{ee} +2\xi G_{ev} -m^2 G_e) \\
        & = -v^7G_e \left[ 2\rho m^2 (G_{ev}+\xi G_{ee}) +\rho^2 (G_{ev}+\xi G_{ee})^2 +m^4 G_{ee}^2 \right]\\
        & \quad +m^2v^7 G_{ee} (G_{ee}G_{vv}-G_{ev}^2) \\
        & = -v^7 G_e \left(   m^2 G_{ee} +\rho (G_{ev} +\xi G_{ee})\right)^2 +m^2v^7 G_{ee} (G_{ee}G_{vv}-G_{ev}^2).
    \end{align*}
    This together with $q_{EE}=-vG_{ee}$ gives
    \begin{align*}
        B  & = q_{EE} A  -v^7 G_e \left(   m^2 G_{ee} +\rho (G_{ev} +\xi G_{ee})\right)^2+m^2v^7 G_{ee} (G_{ee}G_{vv}-G_{ev}^2) \\
        & =v^5 G_eG_{ee}(E^2 G_{ee} +2E G_{ev} +G_{vv}) -v^5 G_e \left( E G_{ee} +G_{ev}\right)^2\\
        & =v^5 G_e (G_{vv}G_{ee}-G_{ev}^2) 
         >0.
    \end{align*}
    We thus conclude that $q$ is convex.
\end{proof}

\section{Internal energy}
Recall that from the Boltzmann kinetic equation with \\$f(t, x, \xi)$ as the probability density of gas particles at position $x$ and moving with velocity $\xi$,
the corresponding energy takes the form
$$
E=\frac{1}{2} \int |\xi|^2 f(t, x, \xi)d\xi,
$$
with
$$
\rho= \int f(t, x, \xi)d\xi, \quad \rho \mathbf{u}= \int \xi f(t, x, \xi)d\xi.
$$
This relation between kinetic description and the averaged fluid variables allows one to show $E \geq \frac{1}{2}\rho |\mathbf{u}|^2$, since
\begin{equation*}
    \frac{1}{2}\rho |\mathbf{u}|^2  = \frac{1}{2}\frac{ \left|\int \xi f(t, x, \xi)d\xi\right|^2}{ \int f(t, x, \xi)d\xi}
     \leq \frac{1}{2} \int |\xi|^2 f(t, x, \xi)d\xi =E.
\end{equation*}
Motivated by this fact, we introduce the following quantity
$$
R=E- \frac{|\mathbf{m}|^2}{2\rho}=\rho e. 
$$
While one can  verify that this function is concave in $\mathbf{w}$, the question here is of course to decide whether $\{\mathbf{w}| R>0, \rho>0\}$ is invariant by the Euler dynamics.

Now we make use of the inequality (\ref{ineq:pve_gamma})
\begin{equation}\label{bb}
    \frac{vG_v}{eG_e} < 2\gamma.
\end{equation}
With this inequality we are able to show that $e>0$, hence $R>0$, is invariant. More precisely we must demand $\gamma$ to be bounded.
But this is true for many EOS and away from vacuum. Formally from both the energy equation and the momentum equation we deduce
$$
\partial_t (\rho e)+\nabla_x \cdot (\rho e \mathbf{u})+P\nabla_x \cdot \mathbf{u}=0.
$$
For smooth solutions we have
$$
\partial_t e +\mathbf{u} \cdot \nabla_x e= \left(-\nabla_x \cdot \mathbf{u} \frac{P}{\rho e} \right)e, 
$$
where the ratio on the right hand side can be written as
\begin{equation*}
    \frac{P}{\rho e}=  \frac{vG_v}{eG_e}<2\gamma.
\end{equation*}
Thus $e>0$ is invariant.

\section{IRP limiter}
For the Euler equations with a general EOS, we have identified a useful invariant region (\ref{ir}), i.e.,
\begin{equation*}
    \Sigma=\{\mathbf{w}\,|\quad \rho>0,\; R>0,\;  q<0\},
\end{equation*}
where
$$
R=E- \frac{|\mathbf{m}|^2}{2\rho}, \quad q=\rho (\inf s_0 - s).
$$
In this section we discuss techniques of how to limit numerical solutions obtained from a high order finite-volume type schemes back to $\Sigma$,
whenever they are out of $\Sigma$ at some solution points.

We first recall the general explicit limiter introduced in \cite{JL18b}. Assume the multi-dimensional system of conservation laws admits an invariant region $\Sigma$, characterized by
\begin{equation*}
\Sigma=\{\mathbf{w}\big | \quad U(\mathbf{w})\leq 0\},
\end{equation*}
where $U$ is convex. Denote the interior of $\Sigma$ by $\Sigma_0$. A key fact we have used is that for any bounded domain $K$, the averaging defined by
$$
\bar{\mathbf{w}}=\frac{1}{|K|}\int _K\mathbf{w}(x)dx
$$ 
is a contraction operator.
\begin{lemma}[ \cite{JL18b}]
    Let $\mathbf{w}(x)$ be non-trivial piecewise continuous vector functions.
    If $\mathbf{w}(x)\in \Sigma$ for all $x\in K \subset \mathbb{R}^d$ and $U$ is strictly convex, then $\bar{\mathbf{w}}\in \Sigma _0$ for any bounded domain $K$.
\end{lemma}
This lemma sets the foundation for using the domain average as a reference to limit the obtained solution polynomials, through a linear convex combination as  in \cite{ZhangShu2010a, ZhangShu2012}.
If we consider a system, the question of particular interest is whether the limiting approximation is still high order accurate.

Let $\mathbf{w}_h(x)$ be a sequence of vector polynomials over $K$, which is a high order accurate approximation to the function $\mathbf{w}(x)\in \Sigma$.
Assume $\bar{\mathbf{w}}_h\in \Sigma _0$, but $\mathbf{w}_h(x)$ is not entirely located in $\Sigma$. We construct
\begin{equation*}
    \tilde{\mathbf{w}}_h(x)=\theta \mathbf{w}_h(x)+(1-\theta)\bar{\mathbf{w}}_h,
\end{equation*}
where $\theta \in (0,1]$ is defined by $\theta =\min\{1,\theta _1\}$, where
\begin{equation*}
    \theta _1=\frac{U(\bar{\mathbf{w}}_h)}{U(\bar{\mathbf{w}}_h)-U^{\max}_h}, \quad U^{\max}_h=\max _{x\in K}U(\mathbf{w}_h(x))>0.
\end{equation*}
If $ \Sigma =\bigcap \limits ^M \limits _{i=1}\{\mathbf{w}\big |\quad U_i(\mathbf{w})\leq 0\}$, then the limiter parameter needs to be modified as
\begin{equation*}
    {\theta =\min \{1,\theta _1, \cdots, \theta _M\}}.
\end{equation*}
This reconstruction has been shown to satisfy three desired properties.
\begin{theorem}\cite{JL18b}\label{BigThm}
    The reconstructed polynomial $\tilde{\mathbf{w}}_h(x)$ satisfies the following three properties:
    \begin{itemize}
        \item[(i)]    the average is preserved, i.e., $\bar{\mathbf{w}}_h=\bar{\tilde{\mathbf{w}}}_h$;
        \item[(ii)]   $\tilde{\mathbf{w}}_h(x)$  lies entirely within invariant region ${\Sigma},  \forall x\in K$;
        \item[(iii)]   order of accuracy is maintained, i.e., if $\| \mathbf{w}_h-\mathbf{w}\|_{\infty}\leq 1$, then
        \begin{equation*}
            \| \tilde{\mathbf{w}}_h-\mathbf{w}\| _{\infty}\leq \frac{C}{|U(\bar{\mathbf{w}}_h)|}\|\mathbf{w}_h-\mathbf{w}\|_{\infty},
        \end{equation*}
        where $C>0$ depends on $\mathbf{w}$ and $\Sigma$.
    \end{itemize}
\end{theorem}
Let $\mathbf{w}^n_h$ be the numerical solution at the $n$-th time step generated from a high order finite-volume-type scheme of an abstract form
\begin{equation*}
    \mathbf{w}^{n+1}_h=\mathcal{L}(\mathbf{w}^n_h), \quad \mathbf{w}^n_h=\mathbf{w}^n_h(x)\in V_h.
\end{equation*}
Provided that the scheme has the following property: there exists $\lambda _0$ and a test set $S$ such that if
\begin{equation*}
{\lambda:=\frac{\Delta t}{\Delta x} \leq \lambda _0 \quad \text{and}\quad \mathbf{w}^n_h(x)\in \Sigma \text{ for }x\in S,}
\end{equation*}
then
$$
{\bar{\mathbf{w}}^{n+1}_h\in \Sigma _0;}
$$
the limiter can then be applied with $K$ replaced by $S_K:S\cap K$, i.e.,
\begin{equation*}
    U^{\max}_h=\max _{x\in S_K}U(\mathbf{w}_h(x)),
\end{equation*}
through the following algorithm:

\textbf{Step 1.}   {\sl Initialization:}  take the piecewise $L^2$ projection of $\mathbf{w}_0$ onto $V_h$, such that
\begin{equation*}
    \langle \mathbf{w}^0_h-\mathbf{w}_0,\phi \rangle =0, \quad \forall \phi \in V_h.
\end{equation*}

\textbf{Step 2.}    {\sl Limiting:} Impose the modified limiter on $\mathbf{w}^n_h$ for $n=0,1,\cdots$ to obtain $\tilde{\mathbf{w}}^n_h$. 

\textbf{Step 3.}  {\sl Update} by the scheme:
\begin{equation*}
    \mathbf{w}^{n+1}_h=\mathcal{L}(\tilde{\mathbf{w}}^n_h).
\end{equation*}
Return to Step 2.

For multi-dimensional Euler equations (\ref{eu}),
we use the invariant region defined in (\ref{ir}) to formulate the IRP limiter which serves as a basis for constructing high order IRP numerical schemes for (\ref{eu}).

A limiter as such was first reported in \cite{JL18a} for one-dimensional Euler equations, and in \cite{JL19} for the isentropic gas dynamics.
The limiter in \cite{JL18a} is explicit and simultaneously preserves the positivity of density and pressure and also a minimum principle for the specific entropy \cite{Ta86}.
\section{An example EOS}
It was shown previously that the polytropic gas EOS admits an invariant region and it easily fulfills the conditions formulated above.
Here we want to present another EOS which fulfills the requirements. To this end we will consider the so-called Tait EOS, cf.\ \cite{Dymond1988}. Let the specific internal energy be given as follows
\begin{align}
   & e = F(s,v) = A(v - v_r) + B\Phi(v) + \frac{1}{2C}\left[(s - s_r) - D(v - v_r)\right]^2\notag\\
    &\hphantom{e = F(s,v)}+ \theta_r\left[(s - s_r) - D(v - v_r) + C\theta_r\right]+ e_r,\label{eq:e_eos}\\
    \Phi(v) &= 
        \begin{cases}
            \ln\dfrac{v_r}{v},\;\nu = 1\\
           \dfrac{1}{1 - \nu}\left(\dfrac{1}{v_r^{\nu - 1}} - \dfrac{1}{v^{\nu - 1}}\right),\;\nu > 1
        \end{cases}\,.
\end{align}
The constants $A,B,C,D$ will be specified later.
The quantities $e_r,v_r, s_r$ and $\theta_r$ denote reference constants of the quantities which may be chosen at a given reference state in the phase space of the material at hand.
The exponent $\nu \geq 1$ can also be chosen according to the specific material under consideration. Now we calculate the pressure and the temperature according to the formulas given above, i.e.,
\begin{align}
    P &= -\frac{\partial F}{\partial v}(s,v) = -\left(A + B\Phi'(v) - \frac{D}{C}\left[(s - s_r) - D(v - v_r)\right] - D\theta_r\right),\label{eq:p_tait}  \\
    \theta &= \frac{\partial F}{\partial s}(s,v) = \frac{1}{C}\left[(s - s_r) - D(v - v_r)\right] + \theta_r.\label{eq:T_tait}
\end{align}
Using equation (\ref{eq:T_tait}) for the temperature we can express the pressure as a function of the volume and the temperature, i.e.,
\begin{equation}
    P = f(v,\theta) = D(\theta - \theta_r) - A - B\Phi'(v).\label{eq:p_tait_vT}
\end{equation}
Such a relation is called thermal EOS. We now choose
\[
 A := K_r - p_r,\,B := K_rv_r^\nu,\quad\text{and}\quad\bar{p}(\theta):= D(\theta - \theta_r) + p_r.
\]
The constant quantity $K_r > 0$ denotes the modulus of compression at a given reference state.
The constants $K_r,v_r$ may be chosen as the saturation values of the material under consideration for a given reference temperature $\theta_r$, cf. \cite{Tf18}.
Accordingly $\bar{p}(\theta)$ can be obtained as the linearization of the saturation curve near a given temperature. A similar approach is used in \cite{HT2015}.
Usually the saturation curve is given as the saturation pressure being a function of the temperature
\[
  p_0 = f_0(\theta),
\]
e.g. \ \cite{Landau1987,Tf18}, and thus we would obtain
\[
  \bar{p}(\theta) = f_0'(\theta_r)(\theta - \theta_r) + p_r.
\]
With these assumptions we can rewrite (\ref{eq:p_tait}) into a more well known form
\begin{equation}
    P = \bar{p}(\theta) + K_r\left[\left(\frac{v_r}{v}\right)^\nu - 1\right].\label{eq:p_tait_vT2}
\end{equation}
The requirement $P > 0$ places a restriction on the allowed region of the state space.
The non-linear Tait EOS is for example used in \cite{Ivings1998,Saurel1999} and may serve as an EOS for liquids.
In order to verify the convexity of the internal energy we calculate the second order derivatives, i.e.,
\begin{align*}
    \frac{\partial^2 F}{\partial v^2}(s,v) = B\Phi''(V) + \frac{D^2}{C},\quad
    \frac{\partial^2 F}{\partial s^2}(s,v) = \frac{1}{C},\quad
    \frac{\partial^2 F}{\partial s\partial v}(s,v) &= -\frac{D}{C}.
\end{align*}
It is easily verified that $F$ is strictly convex if $C > 0$. Indeed $C$ is related to the specific isochoric heat capacity
\[
  c_v = \frac{\partial F}{\partial \theta}(s(v,\theta),v) = C\theta
\]
and can be chosen to be $C = c_{v,r}/\theta_r$. We can also calculate the speed of sound which is given by
\[
  a = \sqrt{v^2\frac{\partial^2F}{\partial v^2}} = \sqrt{v^2\left(B\Phi''(V) + \frac{D^2}{C}\right)}.
\]
For the dimensionless quantities we yield
\begin{align*}
    g &= \frac{Pv}{\theta^2}\frac{\partial^2 F}{\partial s^2} = \frac{Pv}{\theta c_v}, \quad 
    \Gamma = -\frac{v}{\theta}\frac{\partial^2 F}{\partial s\partial v} = \frac{Dv}{c_v},\\
    \gamma &= \frac{v}{P}\frac{\partial^2 F}{\partial v^2} = \frac{v}{P}\left(B\Phi''(V) + \frac{D^2}{C}\right) = \nu\frac{K_rv_r^\nu}{Pv^\nu} + \frac{\Gamma^2}{g},\\
    \mathcal{G} &= -\frac{1}{2}v\dfrac{\dfrac{\partial^3 F}{\partial v^3}}{\dfrac{\partial^2 F}{\partial v^2}} = -\frac{1}{2}v^2\frac{B\Phi'''(V)}{\gamma P}
    = \frac{\nu + 1}{2}\frac{\nu K_r v_r^\nu}{\gamma P v^\nu}.
\end{align*}
Clearly we have $\mathcal{G} > 0$. Thus the isentropes are convex in the $P-v$-plane.
Hence the presented EOS meets all requirements needed above.

\section*{Acknowledgments}
Parts of this research were done during a research stay of Hailiang Liu in Magdeburg.
The authors want to thank Gerald Warnecke for this possibility and his hospitality.
Moreover the authors want to express their gratitude to Marshall Slemrod for helpful discussions.


\medskip
Received February 2021; revised April 2021.
\medskip

\end{document}